\documentclass[12pt]{amsart}
\usepackage{amsfonts}
\usepackage{amsmath}
\usepackage{graphicx}
\usepackage{amssymb}
\newtheorem{theorem}{Theorem}[section]
\theoremstyle{plain}

\newtheorem{defi}{Definition}

\newtheorem{lemma}{Lemma}[section]

\numberwithin{equation}{section}

\newtheorem{prop}{Proposition}[section]
\newtheorem{theo}[prop]{Theorem}

\def\H{\bf H}
\def\A{\bf A}

\def\vphi{\varphi}

\begin{document}
\title[A note on singular time of mean curvature flow]{A note on singular time of  mean curvature flow}
\author{Jingyi CHEN}
\address{Department of Mathematics\\
University of British Columbia\\
Vancouver, B.C., V6T 1Z2\\
Canada}
\email{jychen@math.ubc.ca}
\author{Weiyong He}

\begin{abstract}
We show that mean curvature flow of a compact submanifold in a complete Riemannian manifold cannot form singularity at time infinity if the ambient Riemannian manifold has bounded geometry and satisfies certain curvature and volume growth conditions .

\end{abstract}

\email{whe@math.ubc.ca}
\thanks{The first authors is partially supported by NSERC, and the second author is
partially supported by a PIMS postdoctoral fellowship.}
\date{}

\maketitle

\section{Introduction}

Mean curvature flow develops singularities if the second
fundamental forms of the time dependent immersions become
unbounded. It is well known that mean curvature flow of any closed manifold in the Euclidean space 
develops singularities in finite time. This follows from a maximum principle and barrier argument. In this note, we show,  using integral estimates,  that mean curvature flow cannot form singularity at $t=\infty$ for a class of ambient Riemannian manifolds. More precisely, we prove

\begin{theorem}\label{main}
Let $\Sigma^n$ be a compact manifold and let $M^m$ be a complete
Riemannian manifold with bounded geometry. Suppose that
$F(t):\Sigma\rightarrow M$ satisfies mean curvature flow for
$t\in[0,T)$ and $T$ is the first singular time.  If  $(M, g)$ is Ricci parallel with nonnegative sectional
curvature, and its volume growth satisfies

\begin{equation}\label{E-4-1}Vol(B_p(R))\geq
cR^{m-n+\epsilon}
\end{equation} 
for $R>R_0$, where $\epsilon, c, R_0$ are fixed
positive constants and $B_p(R)$ is the geodesic ball at $p\in M$,  then $T$ has to be
finite. In particular, if $(M, g)$ is analytic, then either
the mean curvature flow $F: \Sigma\rightarrow M$  develops a finite time singularity, or it converges to a compact minimal submanifold in $(M, g)$.
\end{theorem}

A rescaling process is usually applied when
singularities are forming. A sequence of rescaled flows may,
however, move to infinity in ${\mathbb R}^m$ and fail to form a
limit. Particularly, this may happen at type-II singularities if
one scales the flow by the maximum length of the second
fundamental forms at a sequence of times approaching to the first
blow up time. To get compactness, one may consider
the geometric limits for mean curvature flows, in the sense of
Cheeger-Gromov, as Hamilton did for the Ricci flow
\cite{Hamilton95}.

For mean curvature flow, the lower bound on the injectivity radius
follows, unlike the Ricci flow, from the bound of the second
fundamental form ${\A}$. The curvature of a submanifold in the
flow is also bounded if ${\A}$ is bounded, by the Gauss equation.
The smoothness estimate for mean curvature flows \cite{EH} ensures
that all higher derivatives of the second fundamental form are
bounded when ${\A}$ is bounded. These enable one to construct a
limiting mean curvature flow for a sequence of rescaled flows with
a uniform bound on the second fundamental forms. We present a detailed analysis on constructing ancient solutions at any singularity 
and eternal solutions at a type I\!I singularity, although the result (cf. Theorem 2.4) is known, as the results and the arguments will be used in proving Theorem 1.1. 

Combining the geometric limit construction with Hamilton's
monotonicity formula \cite{Hamilton93} for mean curvature flow and Li-Yau's heat kernel estimates \cite{Li-Yau},  we find an upper bound for the first singular time $T$ in terms of volume, and then the volume growth condition imposed on $M$ rules out formation of singularity at infinity.\\\underline{}

\noindent{\bf Acknowledgement:}
Both authors would like to thank Albert Chau for
valuable discussions.

\section{Geometric limit along  mean curvature flow}

Geometric limits of Riemannian manifolds and geometric limits along the Ricci flow are well developed, see  \cite{Hamilton95}, \cite{Morgan-Tian}.
Since we will use the geometric limits  in the formation of singularities along
the mean curvature flow in an essential way, we include some basic facts for completeness.

In this section we do not need to assume the
Riemannian manifolds be compact.

\begin{defi}Let $(M_k, g_k, x_k)$ be a based complete Riemannian
manifold for each positive integer $k$. A geometric limit of the sequence $\{M_k, g_k, x_k\}$ is
a based complete Riemannian manifold $(M_\infty, g_\infty,
x_\infty)$ such that:

(1) there exists an increasing sequence of connected open
subsets $U_k$ of $M_\infty$ exhausting the manifold $M_\infty$,
namely $M_\infty=\cup U_k$ and $U_k$ satisfy
(a) $\overline{U}_k$ is compact,
(b) $\overline{U}_k\subset U_{k+1}$,
(c) $x_\infty\subset U_k$, for all $k$.

(2) for each $k$ there exists a smooth embedding $\vphi_k: (U_k,
x_\infty)\rightarrow (M_k, x_k)$ such that $\varphi_k(x_\infty)=x_k$ and
\[
\lim_{k\rightarrow \infty}\vphi_k^{*}g_k=g_\infty,
\]
where the limit is in the uniform $C^\infty$ topology on compact  subsets of $M_\infty$.

\end{defi}

Similarly, we can define a geometric limit of a sequence of
immersions.
\begin{defi}
Let $F_k: (\Sigma_k, x_k)\rightarrow (N, g, p)$ be a sequence
of immersions with $F_k(x_k)=p\in N$, where $(N, g)$ is a fixed
Rimennian manifold. A geometric limit of the sequence is an
immersion $F_\infty: (\Sigma_\infty, x_\infty)\rightarrow (N,
g, p)$ such that:

(1) there exists an increasing sequence of
connected open subsets $U_k$ of $\Sigma_\infty$, which exhaust the
manifold $\Sigma_\infty$, namely $\Sigma_\infty=\cup U_k$, and satisfy the following for all $k$:
(a)  $\overline{U}_k$ is compact,
(b) $\overline{U}_k\subset U_{k+1}$,
(c) $x_\infty\subset U_k$.

(2) for each $k$ there exists a smooth embedding $\vphi_k: (U_k,
x_\infty)\rightarrow (\Sigma_k, x_k)$ such that $\varphi(x_\infty)=x_k$ and
\[
\lim_{k\rightarrow \infty}F_k\circ \vphi_k =F_\infty,
\]
where the limit is in the uniform $C^\infty$ topology on compact
subsets of $\Sigma_\infty$. In particular,   $ (\Sigma_\infty,
x_\infty, F_\infty^{*}g)$ is a geometric limit of the sequence
$(\Sigma_k, x_k, F_k^{*}g)$ as Rimmennian manifolds.
\end{defi}

A basic fact of geometric limit of Riemmanian manifolds is the
following, which is the $C^\infty$ version of the classical
Cheeger-Gromov compactness.
\begin{theo} \label{Cheeger-Gromov}(Cheeger-Gromov)
Let $(M_k, g_k, x_k)$ be a sequence of connected and based
Riemannian manifolds. Suppose that

(1) for every $R<\infty$, the ball $B(x_k, R)$ has compact
closure in $M_k$ for all $k$ sufficiently  large;

(2) for each integer $l\geq 0$ and each $R<\infty$, there is a
constant $C=C(l, R)$ such that
\[
\left|\nabla^{l}Rm(g_k)\right|\leq C
\]
on $B(x_k,R)$ for all $k$ sufficiently large;

(3) there is a constant $\delta>0$ such that
$inj_{(M_k, g_k)}(x_k)\geq \delta$ for all $k$ sufficiently large.

Then after passing to a subsequence there is a geometric limit
$\{M_\infty, g_\infty, x_\infty\}$ which is a complete Riemmanian
manifold.
\end{theo}

The proof of the above theorem is quite standard (cf. \cite{Morgan-Tian}). For
 evolution equations such as the Ricci flow or the mean
curvature flow, estimates on the higher derivatives of the curvature are the
consequence of the curvature bound, by the smooth estimate. So the key assumption is that
the curvature bound and the injective radius bound. For an
immersion, however, a lower bound on injectivity radius follows from  an upper bound on the second
fundamental form ${\A}$.

A Riemannian manifold $(M,g)$ has bounded geometry if the injectivity radius, the curvatures and the derivatives of the curvatures are uniformly bounded.

\begin{prop}\label{L-2-5}Suppose that
\[
F: \Sigma\rightarrow (M, g)
\] is an immersion where the ambient space $(M, g)$ is a fixed smooth Riemannian manifold with bounded geometry. Suppose that for each $l\geq 0$, there
exists a constant $C=C(l)$ such that
$
\left|\nabla^{l} {\A}\right|\leq C,
$
where $\nabla$ is the covariant derivative of $(\Sigma, F^{*}g)$.
Then the injectivity radius of $(\Sigma, F^{*}g)$ is uniformly bounded from below by a positive constant.
\end{prop}
\begin{proof} We argue by contradiction. Suppose that there exist a sequence of
immersions
\[F_i: (\Sigma_i, x_i)\rightarrow (M, g, F_i(x_i))\] with second
fundamental forms ${\A}_i$ and all their higher derivatives bounded by constants
independent of $i$, but the injectivity radius $\iota_i$ at $x_i\in
(\Sigma_i, F_i^{*}g)$ goes to zero. Consider the sequence
\[
F_i: (\Sigma_i, x_i)\rightarrow \left(M, \frac{1}{\iota^2_i}g,
F_i(x_i)\right).
\]
Then $\left( \Sigma_i, F_i^{*}\left(\iota_i^{-2}g\right),
x_i\right)$ is a sequence of Riemmannian manifolds with bounded
curvature and all higher derivatives of the curvature are also
bounded, by Gauss equation for submanifolds and that
$(M,\iota_i^{-2}g)$ has bounded geometry and $|\nabla^l{\A}_i|\leq
C(l)$. The injective radius of $(\Sigma_i,F^*_i(\iota^{-2}_i g))$
at $x_i$ is $1$. Hence  $\left(\Sigma_i, g_i,x_i\right)$, where
$g_i= F_i^{*}(\{\iota_i^{-2}g)$, converges  in $C^{\infty}$
topology in the sense of Cheeger-Gromov to a geometric limit
$\{\Sigma_\infty, g_\infty, x_\infty\}$, by Theorem
\ref{Cheeger-Gromov}. So there exists an exhausting sequence of
relatively compact open subsets $U_i$ of $\Sigma_\infty$  and a
sequence of $C^{\infty}$ embeddings $\vphi_i$ such that $
\vphi_i^{*}g_i\rightarrow g_\infty$ in $C^{\infty}$ topology on
every compact subset of $\Sigma_\infty$. In particular the
injectivity radius at $x_\infty$ is equal to $1$. To see this,
note that the injectivity radius is equal to the minimum of the
conjugate radius and half of the shortest geodesic loop. In our
case, since the curvature of $g_i$ goes to zero when $i\rightarrow
\infty$, the conjugate radius goes to infinity. Hence, there is a
geodesic loop $l_i$ through $x_i$ in $(\Sigma_i,g_i)$ with length
$2$.  Then $\vphi_i^{-1}(l_i)$ is a sequence of loops through
$x_\infty$ with length converging to $2$. It follows that the
injective radius at $x_\infty$ is less than or equal to 1. It is
clear that the injectivity radius at $x_\infty$ cannot be strictly
less than $1$ as the injectivity radius of $(\Sigma_i,g_i)$ at
$x_i$ is 1.

Note that $ \left(M, \iota^{-2}_ig, F_i(x_i)\right)$ converges to the standard Euclidean space
$\left(\mathbb{R}^m, dx^2, 0\right)$ in $C^{\infty}$
topology on every compact subset. Namely, there exists an
exhausting relatively compact open subsets $V_i$ of $\mathbb{R}^m$
and a sequence of $C^{\infty}$ embeddings $\phi_i$ such that
$\phi_i^{*}\left(\iota_i^{-2}g\right)$ converges to $dx^2$ on every compact subset of $\mathbb{R}^m$ with
$\phi_i(0)=F_i(x_i)$.
Consider  the immersions
\[
\tilde F_i=\phi_i^{-1}\circ F_i\circ \vphi_i: (\Sigma_\infty,
x_\infty)\rightarrow (\mathbb{R}^m, 0).
\]
The second fundamental forms $\tilde{{\A}}_i$ of $\tilde F_i$ are
uniformly bounded and all their higher derivatives are bounded as well
(actually all go to zero), independent of $i$. Hence $\tilde{F}_i$
converges in $C^{\infty}$ topology on compact sets, as a geometric limit, to an immersion
\[
F_\infty: (\Sigma_\infty, x_\infty)\rightarrow (\mathbb{R}^m, 0).
\]
In particular, we have $F_\infty^{*}dx^2=g_\infty$. This statement
is known, however, we include a proof here for completeness. First note that the
injectivity radius at $x_\infty$ is $1$, consider the geodesic ball
$B_1(x_\infty)$ in $(\Sigma_\infty,g_\infty)$, which we can identify with the standard Euclidean
ball $B_1(0)\subset{\mathbb R}^n$ through the exponential map $\exp_{x_\infty}.$
Consider the sequence of immersions
\[
\tilde F_i\circ \exp_{x_\infty}: B_1(0)\rightarrow
\mathbb{R}^m,~~~\tilde F_i(0)=0.
\]
We know that the second fundamental forms of the immersions $\tilde F_i\circ
\exp_{x_\infty}$ are uniformly bounded, this means that the Hessian of the mappings $\tilde F_i\circ
\exp_{x_\infty}$ are uniformly bounded. Also all the higher
derivatives of the second fundamental forms, therefore of the mappings $\tilde F_i\circ\exp_{x_\infty}$,  are uniformly bounded.
It follows that $\tilde F_i\circ \exp_{x_\infty}$ converges in
$B_{1/2}(0)$ to a smooth map $\tilde F_\infty:
B_{1/2}(0)\rightarrow \mathbb{R}^m$ by
Arezella-Ascoli Theorem. We can construct
\[
F_\infty=\tilde F_\infty\circ \exp^{-1}_{x_\infty}
\] in $B_{1/2}(x_\infty). $ To show that we have a limit map
$F_\infty$ on whole manifold $\Sigma_\infty$, we use the
geodesic balls to cover the manifold $\Sigma_\infty$. Note that
for any $y\in \Sigma_\infty$, the injectivity radius is bounded from
below by $d(x_\infty, y)$ since the curvature is uniformly
bounded. The argument then follows from the standard argument of
geometric limit of Riemmanian manifolds by diagonal process. The
reader can refer to \cite{Hamilton95}, \cite{Morgan-Tian} for full
details of the argument.

The second fundamental form $\tilde{\A}_\infty$ of the complete submanifold $F_\infty(\Sigma_\infty)$ is zero because
$$
|\tilde{\A}_i|^2_{g_i}=\iota^{-2}_i|{\A}_i|^2_{F^*_ig}\rightarrow 0,~~\mbox{as $i\to\infty$.}
$$
This implies that $(\Sigma_\infty, x_\infty)$ is a smoothly immersed totally geodesic
submanifold of $\mathbb{R}^m$, hence $(\Sigma_\infty, g_\infty)$
has to be Euclidean itself. But this contradicts
that injectivity radius of $x_\infty$ is equal to $1$.
\end{proof}

By the smoothness property of the mean curvature flow \cite{EH} (the proof holds for general codimension) and Proposition
\ref{L-2-5}, we can get a compactness property along the mean
curvature flow with bounded second fundamental form, similar to
the result of Hamilton \cite{Hamilton95} in the Ricci flow.

\begin{theo}\label{geolimit:flow}
Fix $-\infty\leq T^{'}\leq 0\leq T\leq \infty$ with
$T^{'}<T$. Let $\{\Sigma_k, F_k, x_k\}$ be a sequence of based
mean curvature flows with
\[
F_k(t): \Sigma_k\rightarrow \mathbb{R}^m,~~~F_k{(x_k, 0)}=0.
\]
Suppose that  the lengths of the second fundamental forms
${\A}_k$ of $F_k$ are uniformly bounded above by a constant $C$ independent
of $k$ and time $t$. Then there exists a subsequence of $\{\Sigma_k,F_k,x_k\}$ which converges to
a mean curvature flow $\{\Sigma_\infty, F_\infty(t),  (x_\infty, 0)\}$ as a geometric limit,
where
\[
F_\infty(t): \Sigma_\infty\rightarrow \mathbb{R}^m,
~~~F_\infty(x_\infty, 0)=0,~~~t\in (T^{'}, T),
\]
and $(\Sigma_\infty,F^*_\infty(t)dx^2)$ is a complete Riemannian manifold. 
\end{theo}
\begin{proof}The proof essentially follows Hamilton's argument
\cite{Hamilton95} for the Ricci flow. By Proposition \ref{L-2-5},
the injective radius of $(\Sigma_k, F_k(t)^{*}(dx^2))$ at any
point is uniformly bounded, independent of $t$ and $k$.   Along
the mean curvature flow, the smoothness estimate holds, hence all
the higher derivatives of the second fundamental forms ${\A}_k$ of
$F_k$ are uniformly bounded because $|{\A}_k|$ are uniformly
bounded above. $F_k$ satisfies the mean curvature flow equation,
it follows that all the derivatives of $F_k$ with respect to time
$t$ are also uniformly bounded. Consider the Riemannian manifolds
$(\Sigma_k, F_k(0)^{*}(dx^2), x_k)$. By the assumption, this
sequence has uniformly bounded injective radius and uniformly
bounded curvature and their higher derivatives. It follows that it
sub-converges to a complete Riemannian manifold $(\Sigma_\infty,
g_\infty, x_\infty)$. For any fixed time $T^{'}\leq t_1<0<t_2\leq
T$ and a fixed constant $R$, take a geodesic ball
$B_{2R}(x_\infty)\subset (\Sigma_\infty, g_\infty, x_\infty)$. For
$k$ sufficient large, we can find an embedding
$$
\phi_{k}: B_{2R}(x_\infty)\rightarrow \Sigma_k
$$
 such that
$\phi_k(x_\infty)=x_k$. Define
$$
\tilde F^R_k(t)=F_k(t)\circ\phi_k: B_{2R}(x_\infty)\rightarrow \mathbb{R}^m.
$$
Note $\phi_i$ is time independent. Consider the sequence of immersions
$$
\tilde F^R_k: B_{2R}(x_\infty)\times [t_1, t_2]\rightarrow \mathbb{R}^m
$$
with $\tilde F^R_k(x_\infty, 0)=0$. For simplicity, we can assume
that $2R$ is less than the injective radius and then by using the
exponential map, we can identify $B_{2R}(x_\infty)$ with the
Euclidean ball, as we did in Proposition \ref{L-2-5}. It follows
that all derivatives of $\tilde F_k^R$, as an mapping
from $B_{2R}(x_\infty)\times [t_1, t_2]$ to $\mathbb{R}^m$ with
$\tilde F^R_k(x_\infty, 0)=0$, are bounded. By the classical
Asscoli theorem, it sub-converges to a smooth mapping
$$
F^R_\infty: B_{R}(x_\infty)\times [t_1, t_2]\rightarrow \mathbb{R}^m.
$$
 If $2R$ is larger
 than the injective radius, one applies the covering argument as in \cite{Hamilton95} to show the
convergence in the geodesic ball $B_{R}(x_\infty)$. Now we let $R\rightarrow \infty,
t_1\rightarrow T^{'}, t_2\rightarrow T$ and apply a standard diagonal
sequence argument to obtain a limiting immersion
$$
F_\infty:\Sigma_\infty\times (T^{'}, T)\rightarrow \mathbb{R}^m.
$$
It is clear that $F_\infty (x_\infty, 0)=0$ and $F_\infty$
satisfies the mean curvature flow equation.
\end{proof}

Let $F(t): \Sigma^n \rightarrow (M, g)$ be a smooth mean curvature
flow solution of a compact manifold $\Sigma$ in a smooth
Riemannian manifold $(M, g)$.  We can  use the results above to form a geometric
limit along the mean curvature flow by re-scaling process.
\begin{theo}\label{T-3-1} Let $(M, g)$ be a complete Riemannian manifold with
bounded geometry. If $T$ is the first singular time of the mean
curvature flow $F(t): \Sigma\rightarrow (M, g)$. Then there exists
$(x_i, t_i)$ and $A_i\rightarrow \infty$ such that
\[
F_i(x, s)=F\left(x, \frac{s}{A_i^2}+t_i\right)\rightarrow (M,
A_i^2g)
\]
is a sequence of mean curvature flow solutions, and it
sub-converges to an ancient  mean curvature flow solution for $s\in
(-\infty, 0]$, $
F_\infty(s): \Sigma_\infty\rightarrow \mathbb{R}^m
$
with $|A_\infty(x, s)|\leq |A_\infty(x_\infty, 0)|=1,
F_\infty(x_\infty, 0)=0.$ If the singularity is of type I\!I,
$F_\infty(s)$ can be constructed as an eternal solution. In particular, when $(M, g)=(\mathbb{R}^n, dx^2)$, $\Sigma_\infty(t)$ has at most Eucliden volume growth. 
\end{theo}

\begin{proof}
Suppose that $T$ is the first singular time. For $t<T$, denote
\[A(t)=\max_{p} |{\A}(p, t)|.\] There exist $(x_i, t_i)$ such that
$t_i\rightarrow T$ and
\[
A_i=\max_{t\leq t_i} A(t)= |{\A}(x_i, t_i)|\rightarrow \infty. \]
Consider the sequence of flows defined by
\[
F_i(x, s)=F\left(x, \frac{s}{A_i^2}+t_i\right): \Sigma\rightarrow
\left(M, A_i^2g\right)
\]
for $(x, s)\in \Sigma\times [-A_i^2t_i, 0]$. It is clear that
$F_i(s)$ is still a mean curvature flow solution with
$|A_i(s)|\leq 1$. Set the marked points to be
$
q_i=F\left(x_i,\frac{s}{A_i^2}+t_i\right).
$
It is clear that $\left(M, A_i^2g, q_i\right)$ sub-converges to
the standard Euclidean space $(\mathbb{R}^m, dx^2, 0)$ when
$i\rightarrow \infty$ since $(M, g)$ has bounded geometry.

At $s=0$, the sequence of Riemannian manifolds
 $(\Sigma, F_i(0)^{*}(A_i^2g), x_i)$ sub-converges
smoothly to a Riemannian manifold $(\Sigma_\infty, g_\infty,
x_\infty)$ in the Cheeger-Gromov sense. So there exists a sequence of
exhausting relatively compact open subsets $U_i$ of
$\Sigma_\infty$ and there is a sequence of $C^\infty$ embedding
$$
\vphi_i: U_i\rightarrow (\Sigma, F_i(0)^{*}(A_i^2g)),~~\vphi_i(x_\infty)=x_i.
$$
There also exists a sequence of
exhausting relatively compact open subsets $V_i$ of $\mathbb{R}^m$
and a sequence of $C^{\infty}$ embeddings
$$
\phi_i: V_i\rightarrow(M, A_i^2g),~~\phi_i(0)=F_i(x_i, 0)
$$
such that $\phi_i^{*}\left(A_i^2g\right)$ converges
to $dx^2$ on every compact subset of $\mathbb{R}^m$.

For any fixed $s_0\in (-\infty, 0)$, we
can take $V_i$ such that for any $s\in [s_0, 0]$,
$\phi_i^{-1}\circ F_i(s)\circ \vphi_i(U_i)\subset V_i$, and then we define
$$
\tilde F^R_i(x, s)=\phi_i^{-1}\circ F_i(x, s)\circ \vphi_i: U_i
\rightarrow V_i.
$$
For any fix $R$, by taking $i$ sufficiently large we may assume that $U_i$ contains the
geodesic ball $B_{2R}(x_\infty)$ in $(\Sigma_\infty, g_\infty,
x_\infty)$. It is clear that
\[\tilde F^R_i:
B_{2R}(x_\infty)\subset (\Sigma_\infty, g_\infty,
x_\infty)\rightarrow (\mathbb{R}^m, \phi_i^{*}(A_i^2g),0)
\]
is a mean curvature flow solution with bounded second fundamental form.
By the smoothness property of mean curvature flow, all the higher derivatives of $\tilde
F_k$ and the derivatives with respect to time are also bounded.
Note that the chosen subsequence $\phi_i^{*}(A_i^2g)$ converges to the Euclidean metric $dx^2$
on $\mathbb{R}^m$ when $i$ goes to infinity. Hence the ambient metrics are all
equivalent on any fixed compact subset, especially on $B_R(x_\infty)$. Then we apply the
argument in Theorem \ref{geolimit:flow} to conclude that
the sequence $\tilde F_i$  sub-converges to an immersion
\[
F^R_\infty: B_R(x_\infty)\times [s_0, 0]\rightarrow (\mathbb{R}^m,
dx^2)
\]
with bounded second fundamental form and its higher derivatives in
$B_R(x_\infty)\times [s_0, 0]$.

Now letting $R\rightarrow \infty,$ we obtain a limiting immersion
\[
F_\infty: \Sigma_\infty\times [s_0, 0]\rightarrow (\mathbb{R}^m,
dx^2).
\]
Since $\phi_i$ converges to an isometric embedding independent of
time, it is clear $F_\infty(s)$ still satisfies the mean curvature
flow equation. Taking $s_0\rightarrow -\infty$, we can get an
ancient solution of mean curvature flow
\[
F_\infty: \Sigma_\infty\times (-\infty, 0]\rightarrow
\mathbb{R}^m.
\]
It is clear that $F_\infty(x_\infty, 0)=0,$ and $|A_\infty(x,
s)|\leq |A_\infty(x_\infty, 0)|=1.$

If the mean curvature flow develops a type II singularity at $T$,  one can follow 
Hamilton's work on Ricci flow \cite{Hamilton952} to construct an eternal solution along the mean curvature flow . 

The statement that $\Sigma_\infty(t)$ has at most Euclidean volume growth holds in more general setting. See section 3 for more details of the proof.
\end{proof}

When $n=2$ and ${\H}_\infty\equiv 0$, we have the following
\begin{prop} Let $F(t):\Sigma\to{\mathbb R}^m$ be a smooth mean curvature flow of a compact 2-dimensional surface $\Sigma$ on $[0,T)$.
Let $\Sigma_\infty(s)$ be the geometric limit of $F(t)$ as  in Theorem \ref{T-3-1}.
If ${\H}_\infty\equiv 0$, then $\Sigma_\infty$
has finite total curvature.
\end{prop}
\begin{proof}A complete minimal surface in $\mathbb{R}^m$ for arbitrary $m\geq 3$ has finite total curvature
if and only if it is of finite topological type and has quadratic area growth \cite{Chenqin}.  A surface is of finite
topological type if it has finite genus and finitely many ends.

First,  by Theorem \ref{T-3-1}, the area growth of $\Sigma_\infty$ is at most quadratic.

Second, we show that $\Sigma_\infty$ has finite genus. Since $\Sigma_\infty$ is a geometric limit of $\Sigma$ after
suitable blowing up, there exists a sequence of exhausting open
subsets $U_k$ of $\Sigma_\infty$ and a sequence of embeddings
$\vphi_k$ such that
\[
\vphi_k: U_k\rightarrow \vphi_k(U_k)\subset \Sigma
\]
is a diffeomorphism for each $k$. Then the genus of $U_k$ is less than or equal to that of $\Sigma$. Therefore
$\Sigma_\infty$ has only finite genus since the sequence $\{U_k\}$ exhausts $\Sigma_\infty$.

Finally, we claim that
$\Sigma_\infty$ has only finitely many ends. For any fixed $p\in
\Sigma_\infty$, consider $\Sigma_\infty\backslash B_R(p)$, where
$B_R(p)$ is Euclidean ball in $\mathbb{R}^m$. Let $n_R$ denote the
number of the disjoint components in $\Sigma_\infty\backslash B_R(p)$. If
$\Sigma_\infty$ has infinite many ends, then $n_R\rightarrow\infty$
when $R\rightarrow \infty$. On each component, we can pick up a
point $x_i$ such that the Euclidean distance $d(x_i, p)=2R$ for
$i=1, \cdots , n_R$. We know
\[B_R(x_i)\cap B_R(x_j)=\emptyset\] if $i\ne j$. Now consider $\Sigma_\infty\cap B_{3R}(p)$,
then we know that
\[
\Sigma_\infty\cap \bigcup_{k=0}^{n_R}  B_R(x_k)\subset
\Sigma_\infty\cap B_{3R}(p),
\]
where $x_0=p$. From the monotonicity formula on the area ratio for minimal surfaces in ${\mathbb R}^m$,
\[
area(\Sigma_\infty\cap B_R(x_i))\geq \pi R^2
\]
for $i=0,1,...,n_R$. It follows that
\[
area(\Sigma_\infty\cap  B_{3R}(p))\geq (n_R+1)\pi R^2.
\]
However, we know that $\Sigma_\infty$ has at most Euclidean volume
growth, it follows that
\[
area(\Sigma\cap B_{3R}(p))\leq C(3R)^2.
\]
But this contradicts with $n_R\rightarrow \infty.$ It follows that
$\Sigma_\infty$ has finite many ends.

Therefore, we have shown, by \cite{Chenqin}, that $\Sigma_\infty$ has finite total curvature
$$
\int_{\Sigma_\infty}K = 2\pi l <\infty
$$
where $K$ is the Gauss curvature of $\Sigma_\infty$, $l$ is an
nonnegative integer, and the equality is a classical result of
Osserman on complete minimal surfaces with finite total curvature,
and  $\Sigma_\infty$ is conformally diffeomorphic to a closed
Riemann surface punctured in a finite number of points
\cite{Osserman}.
\end{proof}
\section{Finite time singularity}
Now we assume that $(M, g)$ is Ricci parallel with
non-negative sectional curvatures. Hamilton \cite{Hamilton93}
derived a monotonicity formula for the mean curvature flow with
non-flat ambient space as follows. Note that this coincides to the renowned monotonicity
formula  derived by Husiken when $(M, g)$ is
Euclidean  \cite{Huisken90}. Let $k$ be a solution of the
backward heat equation in $[0, T)$,
\[
\frac{\partial}{\partial t}k=-\triangle_Mk.
\]
Hamilton calculated that
\begin{eqnarray*}
\frac{d}{dt} (T-t)^{(m-n)/2}\int_{\Sigma(t)}
kd\mu=-(T-t)^{(m-n)/2}\int_{\Sigma(t)}\left|\left(H-\frac{Dk}{k}\right)^{\bot}\right|^2kd\mu&&\\
-(T-t)^{(m-n)/2}\int_{\Sigma(t)}g^{\alpha\beta}\left(D_\alpha
D_\beta k-\frac{D_\alpha k D_\beta k
}{k}+\frac{1}{2(T-t)kg_{\alpha\beta}}\right)d\mu.&&
\end{eqnarray*}
If $M$ is Ricci parallel with non-negative positive sectional
curvatures, Hamilton can show that the matrix in the last integral
is non-negative by the Harnack inequality \cite{Hamilton931}. It
follows that
\[
(T-t)^{(m-n)/2}\int_{\Sigma(t)}kd\mu
\]
is non-increasing.  To use the monotonicity formula of Hamilton in an
effective way, we need some properties about the positive
solution of the backward heat equation, which are proved by Li-Yau \cite{Li-Yau} . One can also find the proof in Schoen-Yau \cite{Schoen-Yau}.
\begin{lemma}\label{L-4-1} Let $M$ be a complete manifold with non-negative
Ricci curvature. Let $u$ be a positive solution of the backward
heat equation in $[0, T)$ \[ \frac{\partial }{\partial
t}u=-\triangle _M u
\]
with $\int_M u\equiv 1$ and $u(T)$ is the $\delta$ function
centered at $p$. Then we have

(1).
\[
c (T-t)^{m/2} \leq u(x, t)
\]
when $t\rightarrow T$ and $d(x, p)\leq \sqrt{T-t}$.

(2). For $t>0$ and $x\in M$, we have \[u(x, T-t)\leq C(\delta, m)
Vol_p^{-1/2}(\sqrt{t})Vol_{x}^{-1/2}(\sqrt{t})\exp\left(-\frac{d^2(x,
p)}{(4+\delta)t}\right).\]
\end{lemma}

Now we are in the position to prove Theorem \ref{main}.

\begin{proof} Keep the same notations as in Theorem \ref{T-3-1}. Suppose that $T=\infty$ is the first singular time. This means that the flow exists for all time and as time approaches infinity the second fundamental form becomes unbounded. Denote
$A(t)=\max_{p} |{\A}(p, t)|.$ There exist $(x_i, t_i)$ such that
$t_i\rightarrow\infty$ and $
A_i=\max_{t\leq t_i} A(t)= |{\A}(x_i, t_i)|\rightarrow \infty$ as $i\to\infty$. 
 Denote $p_i=F\left(x_i,
t_i\right).$ By Theorem \ref{T-3-1} we can
construct an ancient solution of mean curvature flow
$$F_\infty(s): \Sigma_\infty\rightarrow \mathbb{R}^m$$
for $s\in (-\infty, 0]$, where $F_\infty(x_\infty, 0)=0$. 

Now we consider the volume growth of $\Sigma_\infty(0)$. For any $R$ fixed, we have
\begin{equation}\label{E-3-1}
\int_{\Sigma_\infty(0)\cap
B_{R}(0)}d\mu_\infty(0)=\lim_{i\rightarrow
\infty}\int_{\Sigma_i(0)\cap B^i_R(p_i)}d\mu_i(0),
\end{equation}
where $B_R(0)$ is the radius $R$ ball in $\mathbb{R}^m$, while
$B^i_R(p_i)$ is the radius $R$ ball in $(M, A_i^2g)$. It
is clear from  the rescaling process that
\begin{equation}\label{E-3-2}
\int_{\Sigma_i(0)\cap B^i_R(p_i)}d\mu_i(0)=\int_{\Sigma(t_i)\cap B^i_{\frac{R}{A_i}}(p_i)}A_i^nd\mu(t_i),
\end{equation}
where $B^i_{\frac{R}{A_i}}(p_i)$ is the radius
$\frac{R}{A_i}$ ball in $(M, g)$. By Lemma
\ref{L-4-1} $(1)$,  we have
\begin{eqnarray}\label{e-4-2}
\int_{\Sigma(t_i)\cap B^i_{\frac{R}{A_i}}(p_i)}A_i^nd\mu(t_i)&=&R^n\int_{\Sigma(t_i)\cap
B^i_{\frac{R}{A_i}}(p_i)}\frac{1}{R^n/A_i^n}d\mu(t_i)\nonumber\\
&=&R^n (T_i-t_i)^{(m-n)/2}\int_{\Sigma(t_i)\cap
B^i_{\frac{R}{A_i}}(p_i)}\frac{1}{(T_i-t_i)^{m/2}}d\mu(t_i)\nonumber\\
&\leq& CR^n (T_i-t_i)^{(m-n)/2}\int_{\Sigma(t)}k_i(x, t_i)d\mu
(t_0),
\end{eqnarray} where
$T_i=\frac{R^2}{A_i^2}+t_i$ and $k_i(x, t)$ for $t\in [0, T_i]$ is the solution of the
backward heat equation
\[
\frac{\partial }{\partial t}k_i=-\triangle_M k_i
\]
with $\int_Mk_i\equiv 1$ and $k_i(x, T_i)$ is the delta function
centered at $F(x_i, t_i)$. Namely, on $(M, g)$ and for each $i$, the function  $k_i(x, T_i-t)$ is the 
fundamental solution of the forward heat equation with singularity at 
the point  $F(x_i,t_i)$ in $M$, for $T_i-t\in [0, \infty)$.  

By Hamilton's monotonicity formula for mean curvature flow, we know that
\begin{equation}\label{E-4-2}
(T_i-t_i)^{(m-n)/2}\int_{\Sigma(t_i)}k_i(x, t_i)d\mu (t_i)\leq
T_i^{(m-n)/2}\int_{\Sigma(0)}k_i(x, 0)d\mu (0).
\end{equation}
Note by Lemma \ref{L-4-1} $(2)$, we have that
\begin{equation}\label{E-4-3}k_{i}(x, 0)\leq
CVol_x^{-1/2}(\sqrt{T_i})Vol_{F(x_i, t_i)}^{-1/2}(\sqrt{T_i}).
\end{equation}

 If $(M, g)$ satisfies the volume growth
condition (\ref{E-4-1}), by (\ref{E-4-2}) and (\ref{E-4-3}), we
get that
\begin{equation}\label{E-3-6}
(T_i-t_i)^{(m-n)/2}\int_{\Sigma(t_i)}k_i(x, t_i)d\mu (t_i)\leq
CVol(\Sigma_0)T_i^{-\epsilon}.
\end{equation}
Since $T_i\rightarrow \infty$ when $i\rightarrow \infty$,  by (\ref{E-3-1}) -- (\ref{E-3-6}), we get that
\begin{equation}\label{E-3-7}
\int_{\Sigma_\infty(0)\cap B_{R}(0)}d\mu_\infty(0)\leq
CT_i^{-\epsilon}Vol(\Sigma_0)R^n.
\end{equation}
Since $\epsilon>0$, the right hand side of (\ref{E-3-7}) goes to zero as $i\to\infty$, but this is impossible. It follows that if $T$ is a singular time it must be finite.  In other words,  if the mean curvature flow solution $F(t)$ exists for
all time, the second fundamental form is uniformly bounded.  In particular, if $(M, g)$ is analytic and there is no finite time singularity,  the mean curvature flow $F(t): \Sigma\rightarrow M$ converges to a compact minimal submanifold along the flow \cite{Simon}.

When $T$ is finite,  one  can bound $k_i(x, 0)$ in (\ref{E-4-3}) in terms of $T$ since $(M, g)$ has bounded geometry.
It follows from (\ref{E-3-1}) -- (\ref{E-4-3}) that  $\Sigma_\infty(0)$ has at most Euclidean volume growth. Note this does not need $(M, g)$ satisfies  volume growth condition (\ref{E-4-1}). 
The proof for any  $\Sigma_\infty(s)$, $s\in (-\infty, 0]$ is similar.
\end{proof}

\end{document}